\documentclass[a4paper,12pt]{article}

\begin{document}

\begin{center}

{\Large {\bf A Criticism on ``A Mathematician's Apology'' by G. H.
Hardy
\\ }}

\end{center}

\vspace{5mm}

\hspace{1in} \parbox[t]{4in}{ {\em ``$\ldots$  modern science
$\ldots$ has recognised the anthropomorphic origin and nature of
human knowledge. $\ldots$ it has recognised that man is the
measure of all things, and that there is no other measure.''} }

\hspace{9cm} \parbox[t]{2in}{{\small {\sc Tobias Dantzig}} }

\vspace{3mm}

Among the greatest scientists in the world, there are many who
consider science, from an idealistic point of view, as an {\em
anonymous} human achievement, hence, something that mankind should
be proud of (contrasting by the way, with many other shames of our
species!). This is the point of view which I want to adopt in this
article: consider science in a classic and pure context,
independent of industrial technology or military motivations.
Following this idea, science is  evolving naturally as a
consequence of the curiosity of the mind on finding out how things
in the world around us work, and how they behave: Is it possible
that there are general laws? Are there some principles governing
the apparent chaos?  Is it possible that the apple which falls is
part of the same natural phenomenon as the Moon which rotates
around the Earth but does not fall?

 The search for such kind of answers is rather exciting for these
 intellects. Two hundred years ago, Laplace explained this exploration
 in his impressive {\em The System of the World}:  ``One of the
 strongest passions is the love of truth, in a man of genius''
 \cite[Book V, Chap. IV]{L}. Real scientists do not compete with each
 other: the challenge is rather to overcome the limitations and
 ignorance of human beings. Inarguably, some scientific contributions
 carry more weight than others. But, once scientists view their
 research as part of an exhilarating scientific voyage, there is no
 room left for dichotomous attitudes, classifying people as winners or
 losers.  And here starts my criticism of the well-known book of G. H. Hardy
 (1877-1947): ``A Mathematician's Apology''.  In this book, published
 for the first time in 1940, he gives his opinion on the
 mathematical world in twenty-nine short chapters. More recent
 editions are easier to find and include a foreword by C. P. Snow
 \cite{H}. Along this article I am going to point out some ideas
 presented in his book which sound to me either controversial,
 having prejudice or could be argued in a more respectful and
 deferential way.

 Right from the start, he apologizes for his own criticism, claiming
 that: ``exposition, criticism, appreciation, is work for second-rate
 minds.'' It is a rather surprising beginning (what does he imply
 about the author of the {\em Exposition du Syst\^eme du Monde}
 \cite{L}, the masterpiece mentioned above?) We all know about very
 talented philosophers, critics, writers, artists, even journalists
 who had been and are still playing a fundamental role in the
 development of science, arts and humanism in general. So, this comment
 sounds very pretentious coming from a mathematician.

  Still in the first chapter he disdains the speech of Alfred
 E. Housman (1859-1936), Kennedy Professor of Latin in the University
 of Cambridge, in his Leslie Stephen lecture on the 9th of May 1933:
 {\em The Name and Nature of Poetry } when, at the very beginning
 \cite[p. 2]{H1} he modestly referred to his previous speech years before
 in the same Senate-House:
\begin{quote}``In these twenty-two years I have improved in some
 respects and deteriorated in others, but I have not so much improved
 as to become a literary critic, nor so much deteriorated as to fancy
 that I have become one.''
\end{quote}
He was reinforcing what he had said in 1911 in the Cambridge
Inaugural Lecture
 {\em The Confines of Criticism, } about literary criticism
 \cite[p. 27]{H2}. Concerning this quotation, Hardy declares: ``$\ldots$
 deplorable that a great scholar and a fine poet should write like
 this.'' I apologise for wishing to express exactly the same words
 about Hardy's declaration.

My disappointment arises specifically from the fact that the book
was
 written by
 such a great mathematician, who not only left many contributions of
 his own, but also, the unique occidental mathematician who was
 considerate enough to recognise the talent of Ramanujan. A quite
 unusual attitude for that time: help and support bringing to light
 exceptionally talented people who come from not so (scientifically)
 prestigious places. For those who are interested in knowing more
 about the relation between Hardy and Ramanujan, see e.g. the Ranganathan's
 book \cite{R}, where one finds details of Ramanujan's meteoric
 and short carrier, and his depressive and unhealthy life.

 Due to the atrocities of the First World War, Hardy had reasons to condemn
 the application of science in military matters, in particular,
 to reprobate the fact that some research on applied mathematics was supporting directly
 those purposes. Needless to say that applied mathematics is much wider than
 those military purposes (by the way Bertrand Russell knew that, and focused his pacifism
 in a more directed way, up to the point of being imprisoned for pacifist acts during the war).
  For some reason related to this, Hardy was
 very proud for being a {\em pure} mathematician
 (``a real mathematician $\ldots$ the purest of the pure'' as
 C. P. Snow described in the Foreword), I would say, almost to the point of treating
 applied mathematics with prejudice. Nevertheless, ironically, contrasting
  with this stereotype, he also
 became famous due to a beautiful result on applied mathematics.
 In 1908 he sent a 2-page letter  to the editor of the
 {\em Science}  \cite{hardy2} with results (concurrently with the German
 physician W. Weinberg) on how proportions of dominant and recessive
 genetic traits propagate in a large mixed population, the well-known
 Hardy-Weinberg law. This result became centrally important in many
 population genetic problems including hemolytic disease, see e.g.,
 among many others introductory textbooks on the subject, Spiess
 \cite{spiess}. Concerning Hardy's outstanding legacy on pure
 mathematics, his contributions are mainly on theory of Diophantine
 analysis, divergent series, Fourier series, Riemann zeta-function and
 the distributions of primes. His greatest collaborators were
 Littlewood and Ramanujan.

 Back to the book which is the focus of my criticism,
 in Chapter seven he describes his ideas on motivation for scientific
 research. He emphasizes some motives which go, in some sense, against those
 ideas I presented in the first two paragraphs. He claims that,
 besides intellectual curiosity, the inspirations come from
 professional pride, ``ambition, desire for reputation, and the
 position, even the power or the money, which it brings''. I
 agree that these latter points represent part of motivation for many among
 us. But definitely they also represent delicate points which,
 when overcharged  can induce acts which hurt (or could be in the
 borderline of) ethic. Regrettably, he emphasizes the last few motives:

\begin{quote}
``So if a mathematician, or a chemist, or even a physiologist,
were to tell me that the driving force in his work had been the
desire to benefit humanity, then I should not believe him (nor
should I think the better of him if I did)''.
\end{quote}

This phrase could sound as an offence for those working on science
for idealism or for those really working for ``love of truth" or
to ``benefit the humanity", and we know that these people, though
the minority, do exist!

These \'elitist or competitive ideas, in the sense that the only
contributions which matter are made by those among the best, is
pejoratively reinforced along the book. They are expressed in
words like: ``mathematical fame $\ldots$ one of the soundest and
steadiest investment'' (Chap. 8) ,``second-rate mind'' or
``$\ldots$ have done something beyond the powers of the vast
majority of men'' (Chap. 6).

Philosophically speaking, I think nowadays it is hard to find
somebody who agrees with his statement in Chapter 27 concerning a
phrase of Hogben: ``The mathematics which can be used `for
ordinary purposes by ordinary men' is negligible.'' Firstly, on
what concerns the \'elitist aspect of this disdainful phrase,
compare it with the Preface of
 {\em
The Mathematics of Great Amateurs }  \cite{C}, where Coolidge says
that, in his book: ``$\ldots$ the number of men included could
easily be doubled or trebled''. Secondly, on what concerns the
presumptuous aspect of Hardy's phrase, contrast it with the modest
and respectful declaration of Laplace, as simple and deep as that
\cite[Book V, Chap.1]{L} :

\begin{quote}
``$\ldots$ l'ing\'{e}nieuse m\'{e}thode d'exprimer tous les
nombres avec dix caract\`{e}res, en leur donnant \`{a} la fois,
une valeur absolue et une valeur de position; id\'{e}e fine et
importante, qui nous para\^{\i}t maintenant si simple, que nous en
sentons \`{a} peine, le m\'{e}rite.'' \footnote{``$\ldots$ the
ingenious method of expressing all numbers by means of ten
symbols, each symbol receiving a value of position as well as an
absolute value; a profound and important idea which appears so
simple to us now that we ignore its true merit.''}
\end{quote}

This quotation of Laplace opens the second chapter of the
remarkable book {\em Number: the Language of Science} by Tobias
Dantzig \cite{D}, author of the epigraph at the top of this
article. His book received many compliments of 20th century top
scholars, including Albert Einstein. Many other authors also
contrast with  Hardy's idea. Besides the already mentioned Laplace
book \cite{L}, it is hard to prevent to mention the excellent work
of some other illustrious authors (I apologise in advance for
omitting so many of them in this short article) like Morris Kline,
in particular his comprehensive survey {\em Mathematical Thought
from Ancient to Modern Times } \cite{K}, Lancelot T. Hogben (who,
in Chapter 27 and 28 Hardy makes clear that he belongs to a
different school from his own), some contemporaries like David
Fowler, Carl Boyer, or Ian Stewart with his vast work involving
updated mathematical objects. It is equally interesting to note
some introductory texts regarding philosophy of science in
general, like (back to the 17th century) the classic ``dialogue''
written by the founder of the modern physics Galileo Galilei
\cite{GG}, some of H. Poincar\'e's works (\cite{P}, e.g.) and the
more philosophical works of Bertrand Russel and Raymond Wilder. As
we said before, many of these books go in an opposite direction of
the ideas presented by Hardy.

 We all know that there are many books on history and
philosophy of mathematics which are partial, elitists or even
tendentious in many aspects. But what really surprises me is the
gap between Hardy's attitude in his life and the ideas expressed
in his book. Naturally, people develop different concepts,
sensitiveness and points of view along their intellectual carrier.
Nevertheless, Hardy's book, compared to others, humanistically
speaking, reflects a very dry, bitter and thorny philosophy.

 Finally, for those who are interested in getting to know more about the
kind of ``irony'' Hardy enjoyed, have a look in the note by A. M.
Vershik \cite{V}. Without wanting to dislodge Hardy's book from
its established status as a statement on mathematical philosophy
by a thoughtful and articulate mathematician, I recommend that in
reading it we ask ourselves whether some of the ideas presented
are presumptuous and scornful to the point of hurting the
development of science and humanism in general.

\vspace{3mm}

\noindent PAULO R. C. RUFFINO

\vspace{3mm}

\noindent This article was written during a visit to the \\
\noindent Mathematical Institute, University of Oxford, \\ 24-29
St. Giles',
Oxford OX1 3LB, UK.\\
 Supported by FAPESP grant no. 00/04591-3.

\noindent Permanent address: \\
\noindent Departamento de
 Matem\'{a}tica, Universidade Estadual de Campinas,\\ 13.081-970 -
 Campinas - SP, Brazil. e-mail: ruffino@ime.unicamp.br

\end{document}